\documentclass[12pt, a4paper]{amsart}

\usepackage[dvipsnames]{color}

\usepackage{verbatim}

\usepackage[english]{babel}

\usepackage{indentfirst}
\usepackage{float}

\usepackage{graphicx}
\usepackage{epsfig}
\usepackage{color}

\usepackage{amssymb, amsmath, amsthm, amscd, ifthen, amssymb}

\usepackage[all]{xy}

\newtheorem{thm}{Theorem}[section]

\newtheorem{Ex}[thm]{Example}

\theoremstyle{definition}

\newtheorem{dfn}[thm]{Definition}

\title[Moduli spaces of planar pentagonal linkages]{Moduli spaces of planar pentagonal linkages:
combinatorial description}

\author{Irina Gorodetskaya}

\begin{document}
\begin{abstract}
Moduli spaces of planar polygonal linkages admit a cell structure which can be realized as a surgery on the permutohedron.  We present a 3D visualization of the result of the surgery for all types of non-degenerate pentagonal linkages.
\end{abstract}

\maketitle

\section{Introduction}

In this paper we visualize configuration spaces of 5-bar planar polygonal linkages~\cite{panina}, also known as planar mechanisms~\cite{fa}, or hinge constructions~\cite{zvon}. Physically a planar polygonal linkage is a collection of rigid bars, joined consecutively by joints in a closed chain. The bars can rotate freely around the joints in the plane.

In the paper~\cite{Stein} Olivier Mermoud and Marcel Steiner have already realized some visualization of moduli spaces of quadrilateral and pentagonal planar linkages using Morse theory. In this article we use a completely different approach and present a new visualization for pentagonal planar linkages that does not only show the topology of the moduli spaces but also reflects the CW-structure introduced in~\cite{panina}.

Firstly, we give some preliminaries.

\begin{dfn} (see~\cite{ziegler}). The \textit{permutohedron} $\Pi_n$ is defined as the convex hull of all points in $\mathbb{R}^n$ that are obtained by permuting the coordinates of the point $(1,2, \ldots ,n)$. The $k$-faces of  $\Pi_n$ correspond to ordered  partitions  of the set $\{1,2,\ldots,n\}$  into $(n-k)$ non-empty parts. The intersection of two (closed) cells is labeled by the coarsest ordered partition of the set $\{1,2,\ldots,n\}$ that refines both partitions that label the two cells.
\end{dfn}

\begin{dfn} (see~\cite{panina}). A \textit{polygonal  $n$-linkage} is a sequence of positive numbers $L=(l_1,\dots ,l_n)$. We also call $L$ \textit{a closed chain} or a \textit{polygon}. We assume that $L$ satisfies the triangle inequality, which guarantees that the below defined moduli space is nonempty.
\end{dfn}

\begin{dfn} (see~\cite{zvon}). A polygonal linkage is called \textit{generic} if it is impossible to put all its vertices
on a straight line. In other words it is impossible to divide the edges into two groups in such a way that the sums of lengths of edges in both groups are equal:
$$\nexists I \subsetneqq \{1 \ldots n\}: \sum_{i \in I} l_{i} = \sum_{i \in \overline{I}} l_{i}.$$
\end{dfn}

In the paper we treat only generic linkages.

\begin{dfn} (see~\cite{panina}). 
 \textit{A configuration} of $L$ in the Euclidean plane $ \mathbb{R}^2$  is a sequence of points $P=(p_1,\dots,p_{n+1}), \ p_i \in \mathbb{R}^2$ with $l_i=|p_i,p_{i+1}|$ and $l_n=|p_n,p_{1}|$.
 \end{dfn}

\begin{dfn} (see~\cite{panina}). The set $M(L)$ of all  configurations modulo orientation preserving
isometries is \textit{the moduli space, or the configuration space
}of the polygonal linkage $L$. In a generic case it is a smooth manifold (see~\cite{fa}).
\end{dfn}

\textbf{Classification of moduli spaces of pentagonal planar linkages.} In~\cite{zvon} Dimitri Zvonkine classifies the configuration spaces of generic 5-bar polygonal linkages. There are 6 types of moduli spaces, Table~\ref{types} shows the representatives of the edge lengths collections and indicates the corresponding moduli space. In the sequel we use the same representatives.

\begin{table}[H]
\centering 
\begin{tabular}{ccr} 
\hline
\hline
pentagon&   &moduli space\\    
\hline
(1, 1, 1, 1, 3)&   &sphere\\
(1, 1, 1, $\varepsilon$, 2)&   &thorus\\
(2, 2, 1, 1, 3)&   &surface of genus 2\\
(1, 1, $\varepsilon$, $\varepsilon$, 1)&   &2 thori\\
(2, 1, 1, 1, 2)&   &surface of genus 3\\
(1, 1, 1, 1, 1)&   &surface of genus 4\\
\hline
\hline                                     
\end{tabular}
 \bigskip
\caption{Types of moduli spaces for 5-bar generic linkages.} 
\label{types}
\end{table}

\textbf{CW-structure on the moduli space of a polygonal linkage.} In the paper~\cite{panina} Gaiane Panina gives an explicit combinatorial description of  $M(L)$  as  of a CW-complex. The entire construction is very much related (but not equal) to the combinatorics of the permutohedron.

\begin{dfn}
A partition of  $L=(l_1,\dots ,l_n)$ is called \textit{admissible } if the total length of any part does not exceed the total length of the rest.

Instead of  partitions of  $L=(l_1,\dots ,l_n)$ we shall speak of partitions of the symbols $\{1,2, \ldots ,n\}$, keeping in mind the lengths $l_i$.
\end{dfn}

\begin{thm}\label{MainThm}
For a generic planar polygonal linkage there exists a structure of a CW-complex $\mathcal{C}(L)$ on its moduli space $M(L)$.
\begin{enumerate}

\item Open $k$-cells $C$ of the complex $\mathcal{C}(L)$ are labeled by cyclically ordered admissible partitions of the set  $\{1,2, \ldots ,n\}$  into $n-k$ non-empty parts.

We are going to mark the label of the cell $C$ as $\lambda(C)$.

\item In particular, the vertices of the complex (that is, cells of dimension 0) are labeled by cyclic orderings of the set $\{1,2, \ldots ,n\}$.

Their number is always $n!$ because the partition $\{1\}\ldots \{n\}$ is admissible, otherwise $(1 \ldots n)$ would not be a physical polygonal mechanism. We identify them with the elements of the symmetric group $S_{n-1}$ by cutting the cyclical ordering at the $n$-th position and omitting $n$.
 
 \item A (closed) cell $C$ belongs to the boundary of some other (closed) cell
    $C'$  iff $\lambda(C)$ is finer than $\lambda(C')$.
    
\end{enumerate}
\end{thm}

In~\cite{panina} this complex is denoted by $CWM^*(L)$.

We obtain a combinatorial structure that is very much related to the combinatorics of the boundary complex of the permutohedron. They have the same incidence relations, but the sets of cells of the boundary complex and of $\mathcal{C}(L)$ are different.  Another difference is in labeling --- the first set is labeled by linearly ordered partitions and the the second by cyclically ordered ones. 

The theorem implies that there exists a natural bijection between the set of vertices of $\mathcal{C}(L)$ and the set of vertices of the permutohedron $\Pi_{n-1}$. Indeed the cells of the permutohedron have the same labeling --- its vertices are labeled by by elements of $S_{n-1}$.  The bijection $\psi$ maps a vertex of $\mathcal{C}(L)$ to a vertex of $\Pi_{n-1}$ that is labeled by the same element of $S_{n-1}$.

$$\psi : Vert(\mathcal{C}(L)) \rightarrow Vert(\Pi_{n-1}).$$

Not only the complex $\mathcal{C}(L)$ admits a $PL$-structure, but it also can be realized as a polyhedron via the following surgery algorithm:

\begin{enumerate}
    \item Start with the complex $\mathcal{C}(L)$ and the boundary complex of the permutohedron $\Pi_{n-1}$. Realize the vertices of $\mathcal{C}(L)$ as the vertices of
$\Pi_{n-1}$ via the above described mapping $\psi$.
    \item Remove some of the facets of $\Pi_{n-1}$ according to the following rule. Every  face $F$ of $\Pi_{n-1}$ is labeled by some ordered partition
    of $\{1, \ldots ,n-1\}$. Add the one-element  set $\{n\}$ at the end of each partition to make it a cyclically ordered partition of $\{1, \ldots ,n\}$.
 If the result is  non-admissible, remove the face $F$ from the  complex.

    This step gives a realization of all the cells of $\mathcal{C}(L)$
 whose label $\lambda$ contains the one-element  set $\{n\}$.

    \item Add some "diagonal" facets according to the following recipe. Take all the cells $C$ of $\mathcal{C}(L)$ such that the part of its label $\lambda(C)$
    containing $n$  has more than one element. If the result is admissible, patch in the facet spanned by the set of vertices $\psi (Vert(C))$ into the complex.
\end{enumerate}

We apply this algorithm in the next sections.

\newpage

\section{The complex $\mathcal{C}(L)$}

\begin{figure}[H]
\centering
     \includegraphics{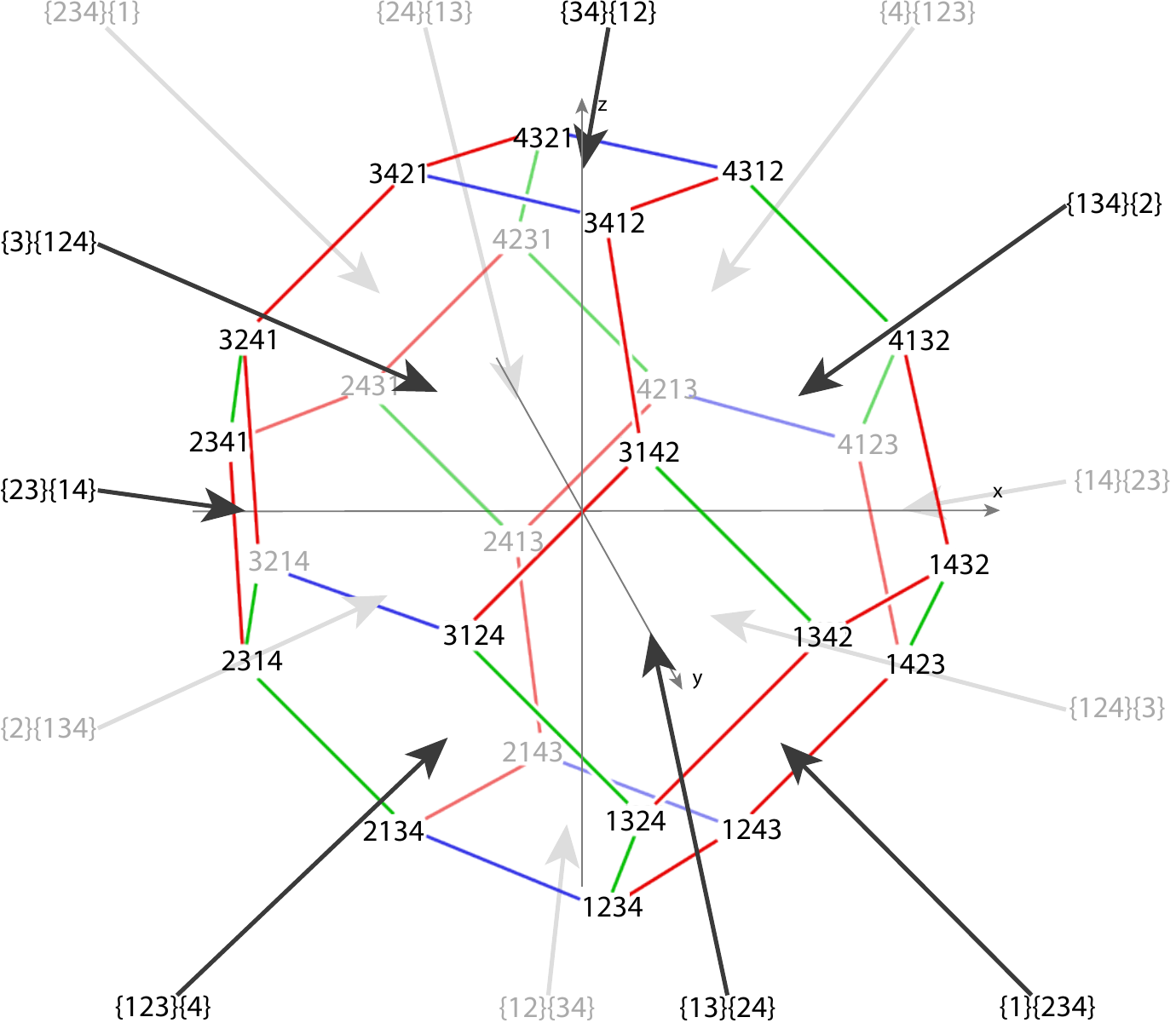}
     \bigskip
\caption{The permutohedron with labeled facets and vertices.}
\label{step_1}
 \end{figure}

Figure~\ref{step_1} depicts a permutohedron $\Pi_{4}$ together with the labels of its facets and vertices.

 \newpage
For each representative from Table~\ref{types} all partitions that correspond to the facets of $\Pi_{4}$ are listed in Table~\ref{step_2}.

\begin{table}[H]
\centering 
\resizebox{\textwidth}{!}{\begin{tabular}{cccccccc} 
\hline
\hline                        
 &Partition&$(1,1,1,1,3)$&$(1,1,1,\varepsilon,2)$&$(2,2,1,1,3)$&$(1,1,\varepsilon,\varepsilon,1)$&$(2,1,1,1,2)$&$(1,1,1,1,1)$\\ 
\hline         
1&$\{1\}\{2,3,4\}\{5\}$&v&v&v&v&v&-\\ 
2&$\{2\}\{1,3,4\}\{5\}$&v&v&v&v&-&-\\  
3&$\{3\}\{1,2,4\}\{5\}$&v&v&-&-&-&-\\
4&$\{4\}\{1,2,3\}\{5\}$&v&-&-&-&-&-\\
5&$\{1,2,3\}\{4\}\{5\}$&v&-&-&-&-&-\\
6&$\{1,2,4\}\{3\}\{5\}$&v&v&-&-&-&-\\
7&$\{1,3,4\}\{2\}\{5\}$&v&v&v&v&-&-\\ 
8&$\{1,2,3\}\{1\}\{5\}$&v&v&v&v&v&-\\ 
9&$\{1,2\}\{3,4\}\{5\}$&v&v&v&-&v&v\\
10&$\{3,4\}\{1,2\}\{5\}$&v&v&v&-&v&v\\
11&$\{1,3\}\{2,4\}\{5\}$&v&v&v&v&v&v\\
12&$\{2,4\}\{1,3\}\{5\}$&v&v&v&v&v&v\\
13&$\{1,4\}\{2,3\}\{5\}$&v&v&v&v&v&v\\
14&$\{2,3\}\{1,4\}\{5\}$&v&v&v&v&v&v\\[1ex]
\hline
\hline      
\\               
\end{tabular}}
\caption{Admissible partitions for step 2 for six representatives from Table~\ref{types}.} 
\label{step_2}
\end{table}

In Table~\ref{step_2} the partitions 1 -- 8 correspond to the hexagonal boundary facets of the permutohedron and the 9 -- 14 correspond to the square ones.

The partitions that consist of the same elements but in different order are simultaneously admissible or non-admissible.

\begin{table}[H]
\centering 
\resizebox{\textwidth}{!}{\begin{tabular}{cccccccc} 
\hline
\hline                        
&Partition&$(1,1,1,1,3)$&$(1,1,1,\varepsilon,2)$&$(2,2,1,1,3)$&$(1,1,\varepsilon,\varepsilon,1)$&$(2,1,1,1,2)$&$(1,1,1,1,1)$\\ 
\hline         
1&$\{3\}\{4\}\{1,2,5\}$ É $\{4\}\{3\}\{1,2,5\}$&-&-&-&-&-&-\\ 
2&$\{2\}\{4\}\{1,3,5\}$ É $\{4\}\{2\}\{1,3,5\}$&-&-&-&-&-&-\\  
3&$\{2\}\{3\}\{1,4,5\}$ É $\{3\}\{2\}\{1,4,5\}$&-&-&-&-&-&-\\
4&$\{1\}\{4\}\{2,3,5\}$ É $\{4\}\{1\}\{2,3,5\}$&-&-&-&-&-&-\\
5&$\{1\}\{3\}\{2,4,5\}$ É $\{3\}\{1\}\{2,4,5\}$&-&-&-&-&-&-\\
6&$\{1\}\{2\}\{3,4,5\}$ É $\{2\}\{1\}\{3,4,5\}$&-&-&-&v&-&-\\
7&$\{3,4\}\{2\}\{1,5\}$ É $\{2\}\{3,4\}\{1,5\}$&-&-&-&-&-&v\\ 
8&$\{2,4\}\{3\}\{1,5\}$ É $\{3\}\{2,4\}\{1,5\}$&-&-&-&-&-&v\\ 
9&$\{2,3\}\{4\}\{1,5\}$ É $\{4\}\{2,3\}\{1,5\}$&-&-&-&-&-&v\\
10&$\{3,4\}\{1\}\{2,5\}$ É $\{1\}\{3,4\}\{2,5\}$&-&-&-&-&v&v\\
11&$\{1,4\}\{3\}\{2,5\}$ É $\{3\}\{1,4\}\{2,5\}$&-&-&-&-&v&v\\
12&$\{1,3\}\{4\}\{2,5\}$ É $\{4\}\{1,3\}\{2,5\}$&-&-&-&-&v&v\\
13&$\{2,4\}\{1\}\{3,5\}$ É $\{1\}\{2,4\}\{3,5\}$&-&-&v&v&v&v\\
14&$\{1,4\}\{2\}\{3,5\}$ É $\{2\}\{1,4\}\{3,5\}$&-&-&v&v&v&v\\
15&$\{1,2\}\{4\}\{3,5\}$ É $\{4\}\{1,2\}\{3,5\}$&-&-&v&-&v&v\\
16&$\{2,3\}\{1\}\{4,5\}$ É $\{1\}\{2,3\}\{4,5\}$&-&v&v&v&v&v\\
17&$\{1,3\}\{2\}\{4,5\}$ É $\{2\}\{1,3\}\{4,5\}$&-&v&v&v&v&v\\
18&$\{1,2\}\{3\}\{4,5\}$ É $\{3\}\{1,2\}\{4,5\}$&-&v&v&-&v&v\\[1ex] 
\hline
\hline      
\\
\end{tabular}}
\caption{Admissible partitions for step 3.} 
\label{step_3}
\end{table}

\newpage

\section{Geometric models}

The figures below (created using Maple) visualize the linkages  $L = (1,1,1,1,1), L = (2,1,1,1,2),L = (2,2,1,1,3), L = (1,1,1,\varepsilon,2)$, and $L = (1,1,1,1,3)$ as polyhedrons. Each figure illustrates stepwise the surgery algorithm.

\begin{Ex} $L = (1,1,1,1,1)$
\begin{figure}[H]
\centering
	\includegraphics{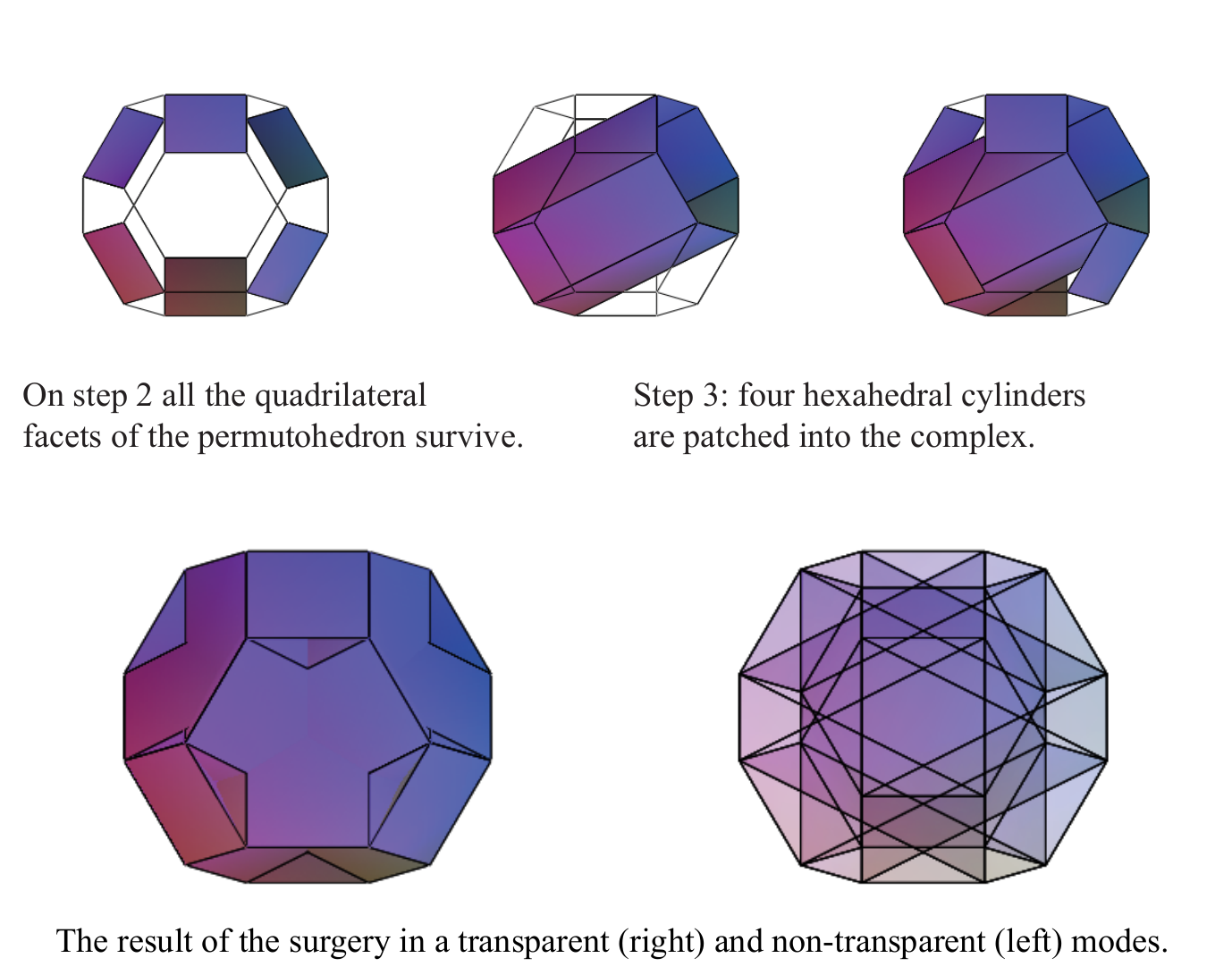}
\caption{$M(1,1,1,1,1)$ is the surface of genus 4.}
\label{genus-4-b}
\end{figure}
\end{Ex}

\newpage

\begin{Ex} $L = (2,1,1,1,2)$
\begin{figure}[H]
\centering
	\includegraphics{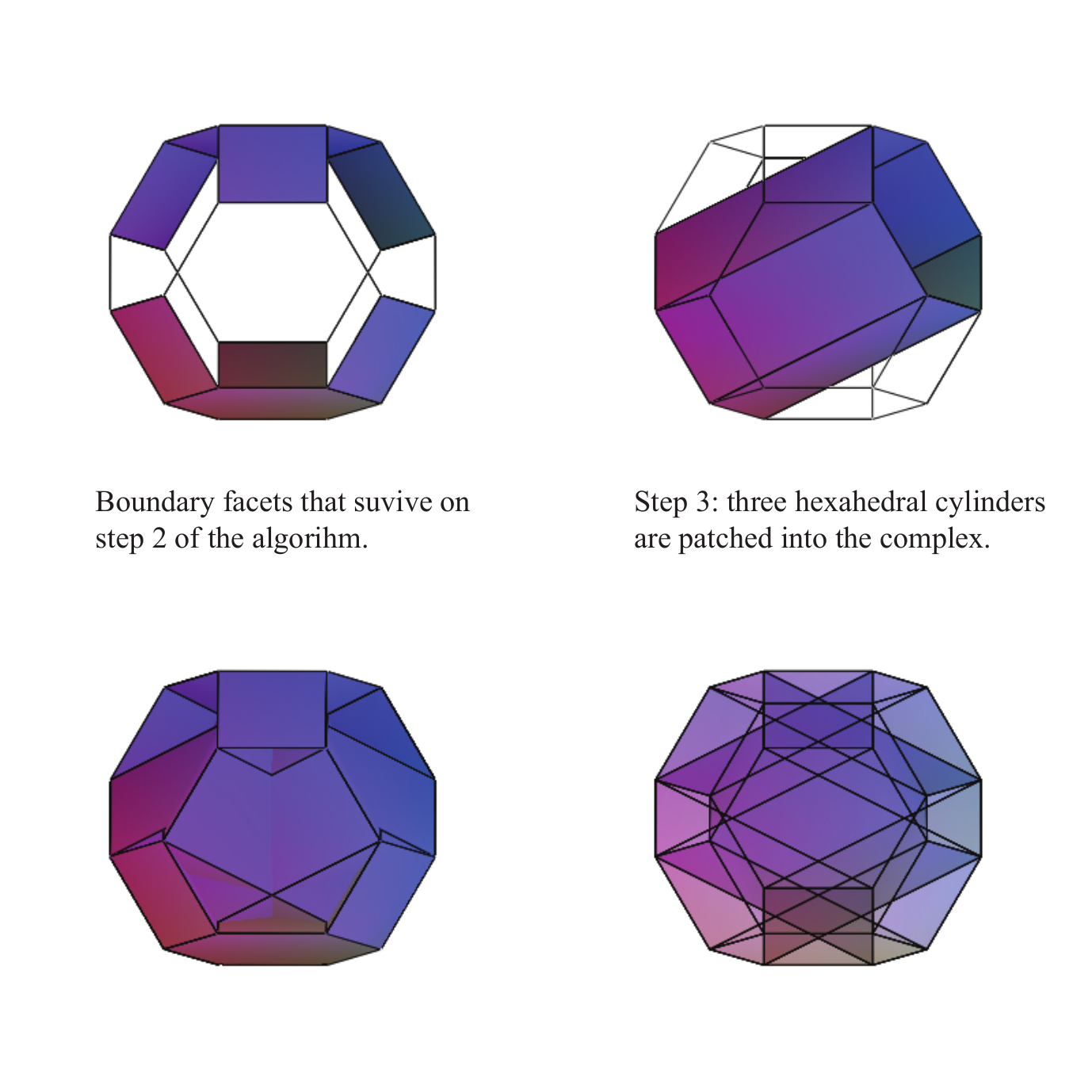}
\caption{$M(2,1,1,1,2)$ is the surface of genus 3.}
\label{genus-3}
\end{figure}
\end{Ex}

\newpage

\begin{Ex} $L = (2,2,1,1,3)$
\begin{figure}[H]
\centering
	\includegraphics{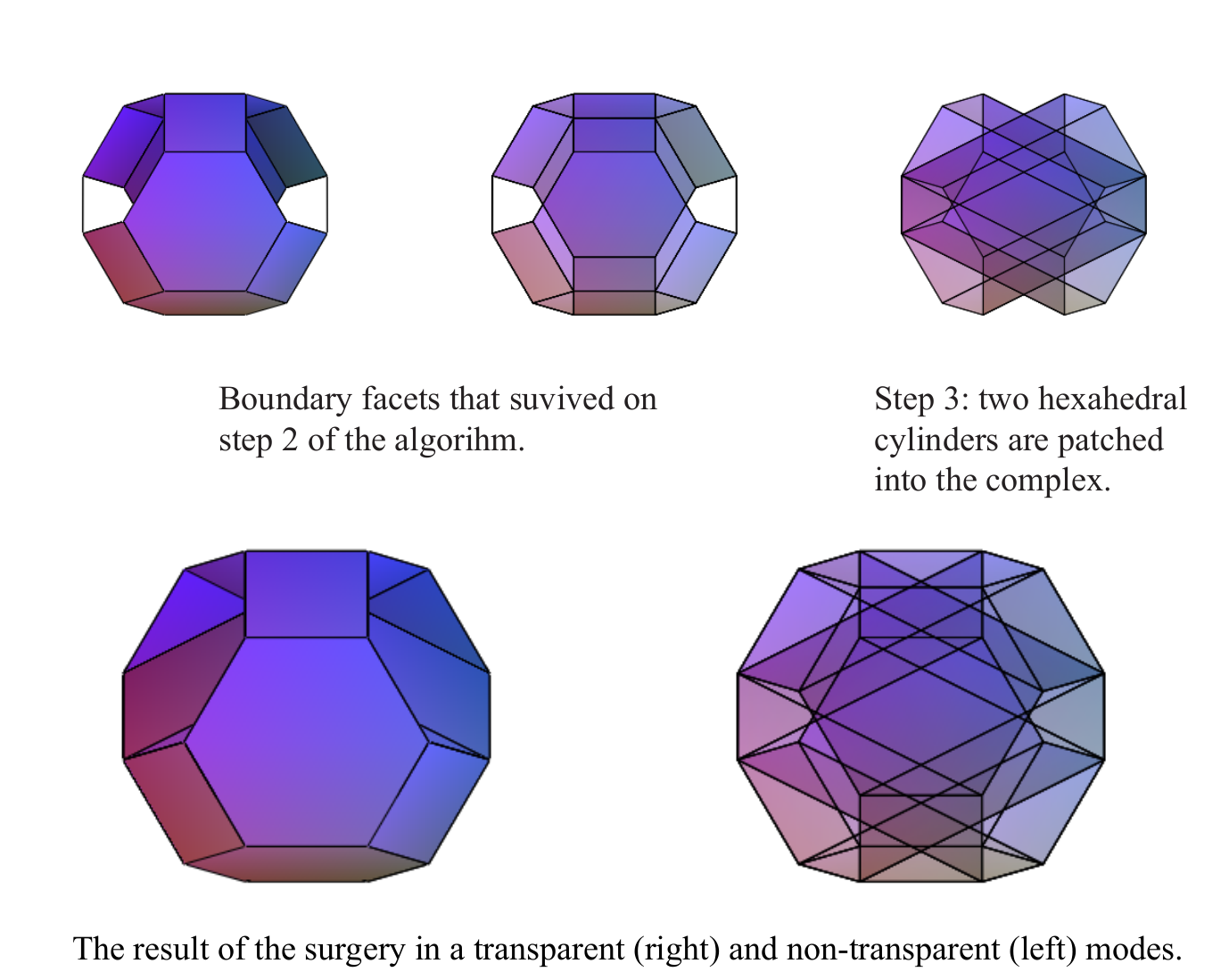}
\caption{$M(2,2,1,1,3)$ is the surface of genus 2.}
\label{genus-2}
\end{figure}
\end{Ex}

\newpage

\begin{Ex} $L = (1,1,1,\varepsilon,2)$
\begin{figure}[H]
\centering
	\includegraphics{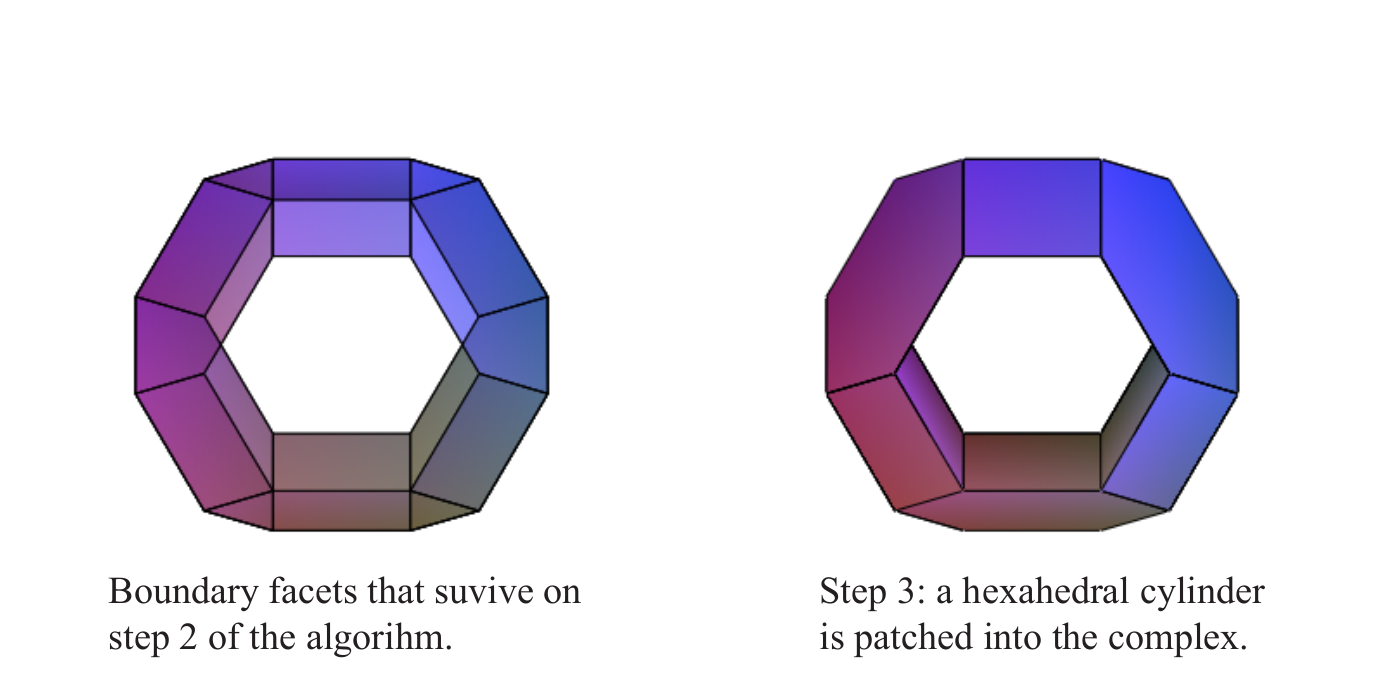}
\caption{$M(1,1,1,\varepsilon,2)$ is the torus.}
\label{thorus}
\end{figure}
\end{Ex}

\begin{Ex} $L = (1,1,1,1,3)$
\begin{figure}[H]
\centering
	\includegraphics{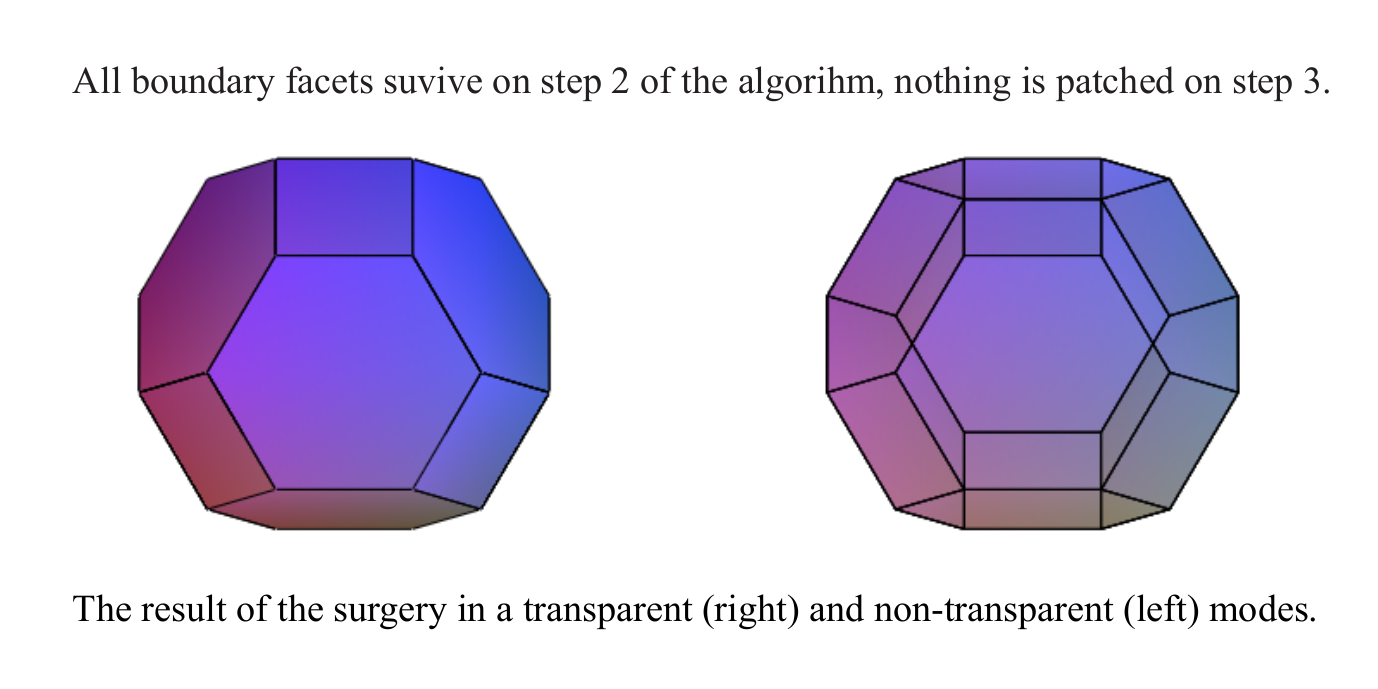}
\caption{$M(1,1,1,1,3)$ is the 2-sphere.}
\label{sphere}
\end{figure}
\end{Ex}

\newpage

The geometric model of the linkage $L = (1,1,\varepsilon,\varepsilon,1)$ is a little bit different from the others. Here on step 2 we not only remove some facets but also some edges. On step 3 we encounter virtual hexagonal facets (see Table~\ref{step_3}). The edges that do not belong to the boundary of any facet are removed. The geometric model breaks up into two symmetric connected components. For better understanding we depict one of them.

\begin{Ex} $L = (1,1,\varepsilon,\varepsilon,1)$
\begin{figure}[H]
\centering
	\includegraphics{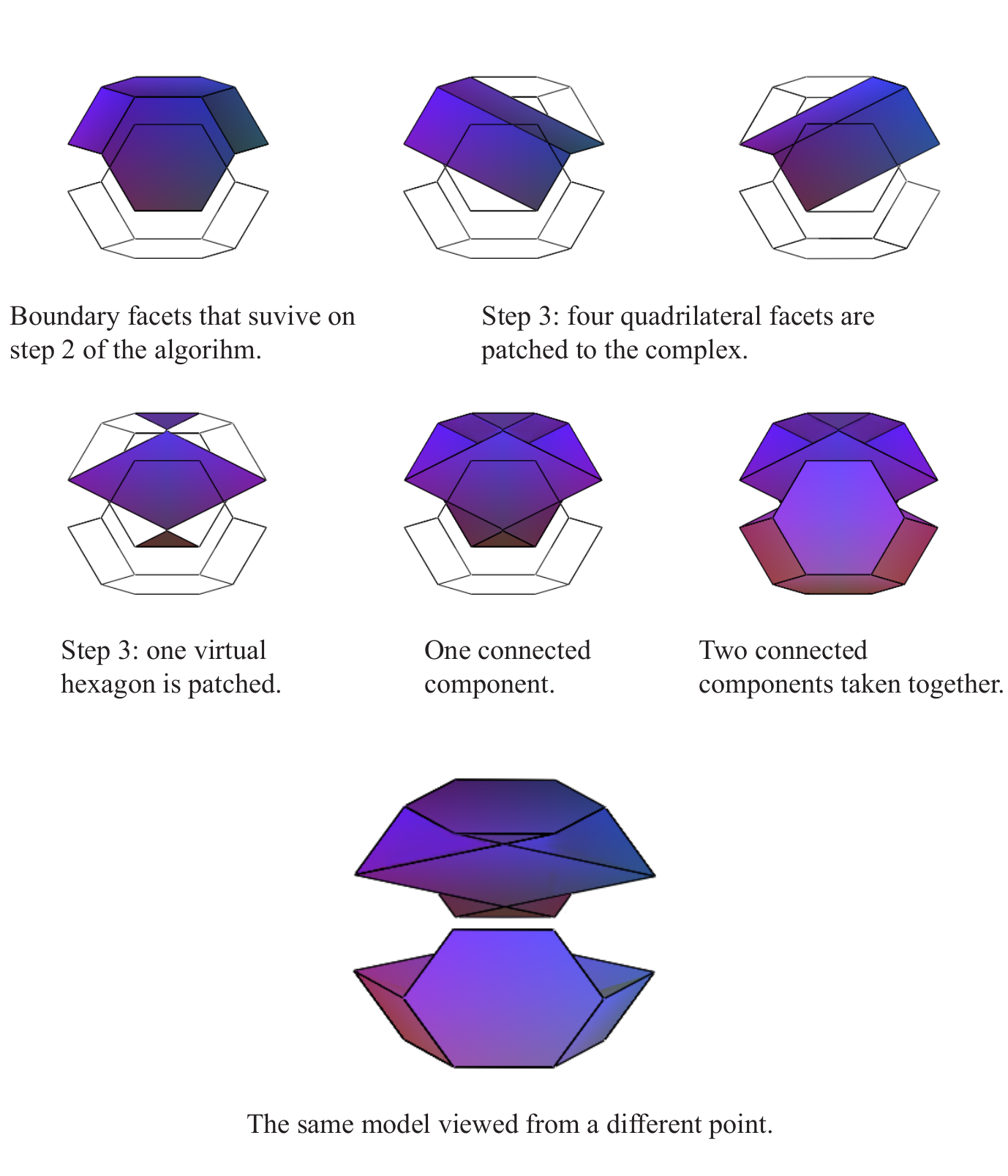}
\caption{$M(1,1,\varepsilon,\varepsilon,1)$ is two tori.}
\label{two-thori}
\end{figure}
\end{Ex}

\end{document}